\documentclass[article]{elsarticle}
\usepackage{lineno,hyperref}
\modulolinenumbers[5]
\journal{--}

\usepackage{amssymb}
\usepackage{amsmath}

\newtheorem{thm}{Theorem}

\newdefinition{rmk}{Remark}
\newproof{pf}{Proof}
\newproof{pot}{Proof of Theorem \ref{thm2}}

\begin{document}

\begin{frontmatter}

\title{On the magnitude of the gaussian integer solutions of the Legendre equation.}
\author{Jos\'e Luis Leal Ruperto}
\address{Departamento de Matem\'atica Aplicada.}
\address{Escuela Polit\'ecnica Superior. Universidad de M\'alaga. Spain.}
\ead{jlleal@uma.es, jllealr@hotmail.com}

\begin{abstract} Holzer \cite{holzer} proves that Legendre's equation 
$$ax^2+by^2+cz^2=0, $$  expressed in its normal form, when
having a nontrivial solution in the integers, has a solution $(x,y,z)$ where 
$|x|\leq\sqrt{|bc|}, \quad |y|\leq\sqrt{|ac|}, \quad |z|\leq\sqrt{|ab|}.$
This paper proves a similar version of the theorem, for Legendre's equation with coefficients $a, b,c$ in Gaussian integers $\mathbb{Z}[i]$ in which there is a solution $(x,y,z)$ where
$$ |z|\leq\sqrt{(1+\sqrt{2})|ab|}.$$
\end{abstract}

\begin{keyword}
\texttt{Legendre Equation; Gaussian Integers} 
\end{keyword}

\end{frontmatter}

\section{Introduction}
Legendre \cite{legendre} proves for the first time that the equation in the integers
\begin{equation}\label{legendre}
ax^2+by^2+cz^2=0 
\end{equation}
expressed in its normal form, that is with $a,b,c$ being square-free and pairwise  coprime, has a nontrivial solution if and only if $a,b,c$ do not all have the same sign and if $-bc$,  $-ca$ y $-ab$ are quadratic residues of $a$, $b$ and $c$ respectively. 
\par
As any multiple of a solution $(x,y,z)$ is a solution, a solution is called primitive if it is simplified and $\gcd(x,y,z)=1$. 
Any equation can be reduced to its normal form: First, if equation 
\begin{equation}\label{legendre1}
a_1x_1^2+b_1y_1^2+c_1z_1^2=0
\end{equation}
is such that
$a_1$, $b_1$, $c_1$ are $a_1=\alpha^2a$,  $b_1=\beta^2b$, $c_1=\gamma^2c$ with $a$, $b$, $c$ being square-free, then there exists a correspondence between the solutions of \ref{legendre1} and those of \ref{legendre}
so that if $(x_1,y_1,z_1)$ is a solution of \ref{legendre1}, then $(\alpha x_1, \beta y_1, \gamma z_1)$ is a solution of \ref{legendre}, and if
 $(x,y,z)$ is a solution of \ref{legendre}, then $(\beta\gamma x, \gamma\alpha y, \alpha\beta z)$ is also a solution of \ref{legendre1}. Second, if equation \ref{legendre} is such that $p$ is a prime divisor of $b$ and $c$, if $(x,y,z)$ is a solution of \ref{legendre}, then $(x,py,pz)$ is a solution of 
 \begin{equation}\label{legendre2}
pax_1^2+\frac{b}{p}y_1^2+\frac{c}{p}z_1^2=0
\end{equation}
and, reciprocally, if $(x_1,y_1,z_1)$ is a solution of \ref{legendre1}, then $(x_1, y_1/p, z_1/p)$ is a solution of \ref{legendre}. As $|pa (b/p)(c/p)|=|abc/p| < |abc|$, after repeating this process as many times as necessary we will come to an equation whose coefficients are pairwise coprime.\par
If $(x,y,z)$ is a primitive solution of \ref{legendre} expressed in its normal form, then we necessarily have that $\gcd(x,y)=\gcd(x,z)=\gcd(y,z)=1$, since in the case $\gcd(x,y)\neq 1$ we would have that $x=px_o$, $y=py_o$ having a common prime factor $p$ with $-cz_o^2=a(px_o)^2+b(py_o)^2$, considering residue classes modulo $p^2$ we would have that $[-cz_o^2]_{p^2}=[p^2]_{p^2}[ax_o^2+by_o^2]_{p^2}=[0]_{p^2}$ which is not possible since $p$ does not divide $z_o$ and $p^2$ does not divide $z_o^2$, nor does $p^2$ divide $c$ which is square-free.
\par
All this is valid in any unique factorization domain. In particular, Gauss integral domain $\mathbb{Z}[i]$ is. 
\par
Samet \cite{samet} proves a similar theorem for equation \ref{legendre} in $\mathbb{Z}[i]$ expressed in its normal form, which has a solution if and only if $bc$, $ca$, $ab$ are quadratic residues of $a$, $b$ and $c$ respectively, omitting the negative sign in the conditions since $i^2=-1$. 
\par
\par
There cannot exist a general formula in terms of coefficients giving a solution to \ref{legendre} as a result. However, on the basis of an already known solution $(x_o, y_o, z_o)$, it is possible to find parametrically the other infinite solutions.
Holzer, in \cite{holzer}, proves that if equation \ref{legendre} in the integers in its normal form has a solution, it is always possible to find one of them $(x,y,z)$ which satisfies 
\begin{equation}
|x|\leq\sqrt{|bc|}, \quad |y|\leq\sqrt{|ac|}, \quad |z|\leq\sqrt{|ab|}
\end{equation}
so that, at worst, if coefficients aren't very large, solutions can be found by trial and error. \par

Nevertheless, no results on the magnitude of the solutions of Legendre's equation in $\mathbb{Z}[i]$ have been found hitherto.\par 
Holzer's theorem, at least in its original form, does not work for Gaussian integers. For instance, equation 
$ix^2+7y^2+z^2=0$ 
has as its smallest solution $(2+2i,1,1)$ where $|x|=|2+2i|=\sqrt{8}>\sqrt{|bc|}=\sqrt{|7|}$. \par
We adapt the proof of Holzer's theorem made years later by Mordell \cite{mordell} in order to prove that at least, the following theorem holds true:

\begin{thm}
The equation  $ax^2+by^2+cz^2=0$ in $\mathbb{Z}[i]$
expressed in its normal form, if it has solutions, has a solution $(x,y,z)$ where
$$ |z|\leq\sqrt{(1+\sqrt{2})|ab|}.$$
\end{thm}

\section{Proof of the theorem}

 The proof will consist of showing that, if there exists a solution $(x_o, y_o, z_o)$ where $(x_o,y_o)=1$ and $|z_o|>\sqrt{(1+\sqrt{2})|ab|}$, then another solution $(x,y,z)$ where $|z|<|z_o|$ can be found on the basis of the first solution.\par
We parametrize the solutions establishing that  
\begin{equation}\label{condicion}
(x_o+tX, y_o+tY, z_o+tZ)
\end{equation}
 is a solution in Gaussian rationals with $X$, $Y$, $Z$ being Gaussian integer parameters and $t$ being a Gaussian rational, $t\neq 0$. By replacing in \ref{legendre}, we get that
\begin{align*}
0=a(x_o+tX)^2+b(y_o+tY)^2+c(z_o+tZ)^2
&=\\
ax_o^2+2ax_otX+at^2X^2+by_o^2+2by_otY+bt^2Y^2+cz_o^2+2cz_otZ+ct^2Z^2
&=\\
ax_o^2+by_o^2+cz_o^2+(aX^2+bY^2+cZ^2)t^2+2t(ax_oX+by_oY+cz_oZ)
&=\\
t((aX^2+bY^2+cZ^2)t+2(ax_oX+by_oY+cz_oZ)).
\end{align*}
As $t\neq 0$, we get that $t=\frac{-2(ax_oX+by_oY+cz_oZ)}{aX^2+bY^2+cZ^2}$. Substituting $t$ in \ref{condicion} we get the solutions of \ref{legendre} in Gaussian rationals;
if we multiply by $aX^2+bY^2+cZ^2$ we get the  integer solutions
\begin{eqnarray}\label{sinsimplificar}
x_o(aX^2+bY^2+cZ^2)-2X(ax_oX+by_oY+cz_oZ)\\
\nonumber
y_o(aX^2+bY^2+cZ^2)-2Y(ax_oX+by_oY+cz_oZ)\\
\nonumber
z_o(aX^2+bY^2+cZ^2)-2Z(ax_oX+by_oY+cz_oZ)
\end{eqnarray}
which can be simplified by using a common divisor $\delta$ to all three expressions in \ref{sinsimplificar}, not necessarily the greatest common divisor, so that we get a solution $(x,y,z)$ given by
\begin{eqnarray}\label{simplificada}
x=\frac{1}{\delta}(x_o(aX^2+bY^2+cZ^2)-2X(ax_oX+by_oY+cz_oZ))\\
\nonumber
y=\frac{1}{\delta}(y_o(aX^2+bY^2+cZ^2)-2Y(ax_oX+by_oY+cz_oZ))\\
\nonumber
z=\frac{1}{\delta}(z_o(aX^2+bY^2+cZ^2)-2Z(ax_oX+by_oY+cz_oZ)).
\end{eqnarray}
We now prove that if 
$$
\delta\,\, \mbox{divides} \,\, c,\,\,\mbox{and also divides}\,\,  Xy_o-Yx_o,
$$
then it divides all three expressions in \ref{sinsimplificar} and therefore $(x,y,z)$ are integers.
To do so, we must previously consider that
 $a$, $b$, $x_o$ and $y_o$ are invertible in the ring $\mathbb{Z}[i]_\delta$. It will be sufficient to notice that $\gcd(\delta, abx_oy_o)=1$. Let's consider the existence of a prime number $p$ being a common divisor of $\delta$ and $abx_oy_o$. As $\delta|c$ and $a$, $b$, $c$ are coprime, then $p$ could only divide $x_oy_o$. Let's consider it divides $x_o$. As $ax_o^2=-by_o^2-cz_o^2$ and $p$ does not divide $b$, it must divide $y_o^2$ and therefore it must also divide $y_o$ which is a contradiction since $\gcd(x_o, y_o)=1$.
\par
 Having said that, as $\delta$ divides $Xy_o-Yx_o$, then $[Xy_o-Yx_o]_\delta=[0]_\delta$. If we isolate $[X]_\delta=[Yx_o/y_o]_\delta $ and replace in $[aX^2+bY^2+cZ^2]_\delta$ and in $[ax_oX+by_oY+cz_oZ]_\delta$
 \smallskip
\begin{align*}
[a(Yx_o/y_o)^2+bY^2+cZ^2]_\delta&=\left[(ax_o^2Y^2+byo^2Y^2+cy_o^2Z^2)/y_o^2)\right]_\delta=\\ 
	=\left[(-cz_o^2Y^2+cy_o^2Z^2)/y_o^2\right]_\delta&=[c]_\delta\left[(-z_o^2Y^2+y_o^2Z^2)/y_o^2\right]=[0]_\delta \\ 
[(ax_o(Yx_o/y_o)+by_oY+cz_oZ]_\delta&=\left[(ax_o^2Y+byo^2Y+cy_o^2Z)/y_o)\right]_\delta=\\ 
								=\left[(-cz_o^2Y+cy_o^2Z)/y_o\right]_\delta &=[c]_\delta\left[(-z_o^2Y+y_o^2Z)/y_o\right]=[0]_\delta \\
\end{align*}
we get multiples of $\delta$ since $c$ also is a multiple of $\delta$. Hence, the expressions in \ref{sinsimplificar} are also multiples of $\delta$, and $(x,y,z)$ is an integer solution.\par
We manipulate equation \ref{simplificada} so as to be able to compare $z$ to $z_o$:
{\small
\begin{align*}
\frac{-\delta z}{cz_o}&=\frac{2Z(ax_oX+by_oY+cz_oZ)}{cz_o}-\frac{z_o(aX^2+bY^2+cZ^2)}{cz_o}\\
				&=Z^2+2Z(\frac{ax_oX+by_oY}{cz_o})^2-\frac{aX^2+bY^2}{c}\\
				&=(Z+\frac{ax_oX+by_oY}{cz_o})^2-(\frac{ax_oX+by_oY}{cz_o})^2-\frac{aX^2+bY^2}{c}\\
				&=(Z+\frac{ax_oX+by_oY}{cz_o})^2-\frac{1}{c^2z_o^2}(aX^2cz_o^2+bY^2cz_o^2+a^2x_o^2X^2+b^2y_o^2Y^2+2abx_oy_oXY)\\
				\intertext{as $cz_o^2=-ax_o^2-by_o^2$,}
				&=(Z+\frac{ax_oX+by_oY}{cz_o})^2+\frac{ab}{c^2z_o^2}(y_o^2X^2-2x_oy_oXY+x_o^2Y^2)\\
				&=(Z+\frac{ax_oX+by_oY}{cz_o})^2+\frac{ab}{c^2z_o^2}(y_oX-x_oY)^2
				\end{align*}
				}
and we get that 
\begin{equation}\label{identidad}
z=z_o\frac{-c}{\delta}\left[(Z+\frac{ax_oX+by_oY}{cz_o})^2+\frac{ab}{ c^2z_o^2}(y_oX-x_oY)^2\right].
\end{equation}
As $gcd(x_o,y_o)=1$, the equation $y_oX-x_oY=\delta$ always has a solution. For our purposes we will choose the values of $X$, $Y$, $Z$ depending on whether $c$ is or is not a multiple of $(1+i)$, as follows:\par
\medskip
 If $c$ is a multiple of $1+i$, we can take $\delta=\frac{1}{1+i}c$, take parameters $X$, $Y$ as any solution of $y_oX-x_oY=\frac{1}{1+i}c$ and the value of $Z$ as the nearest Gaussian integer to the Gaussian rational $-\frac{ax_oX+by_oY}{cz_o}$. That being so, we have that 
 $$\left|Z+\frac{ax_oX+by_oY}{cz_o}\right|<\frac{\sqrt{2}}{2}.$$
 
 This is thus because, geometrically, Gaussian integers form a grid whose edge is equal to 1 and whose diameter is $\sqrt{2}$. Any Gaussian rational located inside a square would be at a distance of half the diameter from a Gaussian integer.\par
 
Let's consider that $|z_o|>\sqrt{(1+\sqrt{2})|ab|}$, taking modulus in \ref{identidad}
 \begin{align*}
|z|&=|z_o||1+i|\left|(Z+\frac{ax_oX+by_oY}{cz_o})^2+\frac{(1+\sqrt{2})ab}{(1+\sqrt{2})(1+i)^2z_o^2}\right|<\\
&<|z_o|{\sqrt{2}}\left((\frac{\sqrt{2}}{2})^2+\frac{1}{(1+\sqrt{2})(\sqrt{2})^2}\right)=|z_o|{\sqrt{2}}\frac{2+\sqrt{2}}{2+2\sqrt{2}}=|z_o|.\\
\end{align*}

  \medskip\par
If $c$ is not a multiple of $1+i$, we choose $X$ and $Y$ as any solution to $y_oX-x_oY=c$. As $\mathbb{Z}[i]_{1+i}=\{[0]_{1+i},[1]_{1+i}\}$, if  $Z$ satisfies the equation
\begin{equation}
[aX+bY+cZ]_{1+i} =[0]_{1+i}
\end{equation}
 then the expressions in \ref{sinsimplificar} are also divisible by $1+i$, as $aX^2+bY^2+cZ^2$ also is, since whatever the value of $\alpha$ in $\mathbb{Z}[i]$ might be, we have that $ [\alpha^2]_{1+i}=[\alpha]_{1+i}$, and then
 \begin{equation*}
 [aX^2+bY^2+cZ^2]_{1+i}=[aX+bY+cZ]_{1+i}=[0]_{1+i}.
\end{equation*}
Furthermore, $2=(1+i)(1-i)$ is also divisible by $1+i$. In this way, it is possible to determine the parity of $Z$, which is either a multiple of $1+i$ or 1 plus a multiple of $i+1$. Let's then choose $Z$ as the Gaussian integer having such parity being nearest to 
$-\frac{ax_oX+by_oY}{cz_o}$ and we will have that 
 $$\left|Z+\frac{ax_oX+by_oY}{cz_o}\right|<1$$
 since the multiples of $1+i$, and those which are not, are located inside a grid whose edge is $\sqrt{2}$ and whose diameter is $((\sqrt{2})^2+(\sqrt{2})^2)^{1/2}=2$.

We replace the divisor $\delta$ in \ref{identidad} by $(1+i)\delta$ and taking modulus,
  \begin{align*}
  |z|&<|z_o|\frac{1}{\sqrt{2}}\left(1+\frac{1}{1+\sqrt{2}}\right)=|z_o|\frac{1}{\sqrt{2}}\frac{2+\sqrt{2}}{1+\sqrt{2}}=|z_o|.\\
  \end{align*}
  By repeating this process, as long as $|z_o|^2>(1+\sqrt{2})|ab|$ we will come to a solution where $|z^2|\leq (1+\sqrt{2})|ab|$.\par
 We must finally show that $z$ can never be $0$, since in such a case the solution would necessarily be the trivial one $(0,0,0)$. If $z$ were 0, that is
 \begin{equation}\label{cero}
z_o(aX^2+bY^2+cZ^2)-2Z(ax_oX+by_oY+cz_oZ)=0
\end{equation}
When isolating, the only possible value of $Z$ would be
\begin{align*}
Z&=\frac{-2(ax_oX+by_oY)\pm \sqrt{4(ax_oX+by_oY)^2+4cz_o(az_oX^2+bz_oY^2)  }}{2cz_o}\\
&=\frac{-(ax_oX+by_oY)}{cz_o}\pm \frac{\sqrt{ab(y_oX-x_oY)^2}}{cz_o}
\end{align*}

However, no integer $Z$ can be obtained in this way, as $\sqrt{ab(y_oX-x_oY)^2}$ is, depending on the case, equal to $\sqrt{ab\frac{c^2}{4}}$ or to $\sqrt{abc^2}$, and with irrational values in both cases, since $a$, $b$ are square-free.
\par
This completes the proof. \par\smallskip
\par

 \bigskip\bigskip\bigskip
\bibliographystyle{elsarticle-num}

\end{document}